\documentstyle[12pt]{article}
\begin{document}
\title{Connections between Hilbert W*-modules \\ and direct integrals}
\author{Michael Frank}
\maketitle

\newtheorem{theorem}{{\sc Theorem}}[section]
\newtheorem{corollary}[theorem]{{\sc Corollary}}
\newtheorem{example}[theorem]{{\sc Example}}
\newtheorem{definition}[theorem]{{\sc Definition}}
\newtheorem{proposition}[theorem]{{\sc Proposition}}
\newtheorem{remark}[theorem]{{\sc Remark}}

\begin{abstract}
Investigating the direct integral decomposition of  von Neumann algebras of
bounded module operators on self-dual Hilbert W*-modules an equivalence
principle is obtained which connects the theory of direct disintegration of
von Neumann algebras on separable Hilbert spaces and the theory of von Neumann
representations on self-dual Hilbert {\bf A}-modules with countably generated
{\bf A}-pre-dual Hilbert {\bf A}-module over commutative separable W*-algebras
{\bf A}. Examples show posibilities and bounds to find more general relations
between these two theories, (cf. R. Schaflitzel's results). As an application
we prove a Weyl--Berg--Murphy type theorem: For each given commutative
W*-algebra {\bf A} with a special approximation property (*) every normal
bounded {\bf A}-linear operator on a self-dual Hilbert {\bf A}-module with
countably generated {\bf A}-pre-dual Hilbert {\bf A}-module is decomposable
into the sum of a diagonalizable normal and of a ''compact'' bounded
{\bf A}-linear operator on that module.
\end{abstract}

The idea to investigate the subject treated in the present paper
arose in discussions with K. Schm\"udgen and  J. Friedrich
at the University of Leipzig. They suggested to the author that
self-dual Hilbert W*-modules over commutative W*-algebras might be closely
connected with direct integrals of measurable fields of Hilbert spaces or,
respectively, with some topologically related objects. Moreover, von Neumann
algebras of bounded module operators on these self-dual Hilbert W*-modules
should be decomposable into direct integrals of measurable fields of von
Neumann algebras in a very easy way.

Following this line appropriate facts have been proved. One gets a new view
on the nowadays well-known theory of
direct integral decomposition of von Neumann algebras {\bf M} on separable
Hilbert spaces. This theory is shown to be equivalent to the theory of von
Neumann representations of W*-algebras {\bf M} on self-dual Hilbert W*-modules
$\cal H$ over W*-subalgebras {\bf B} of the center of {\bf M}, where {\bf B}
has to be separable and the Hilbert {\bf B}-modules have to possess countably
generated {\bf B}-pre-dual Hilbert {\bf B}-modules. The most interesting point
is that the basic structures, Hilbert W*-modules and direct integrals of
Hilbert spaces, are quite different. However, this equivalence
will not be preserved turning to direct integrals of von Neumann
algebras on non-separable Hilbert spaces, in general. It would be interesting
to make further consi\-derations in this direction taking in account recent
results of R. Schaflitzel \cite{Schaflitzel:90/1,Schaflitzel:90/2},
\linebreak[4] P. Richter \cite{Richter:91} and other authors
\cite{Kehlet:78,Marechal:69,Vesterstrom/Wils:70,Wils:70}. Applicating this
equivalence principle, a new result is found generalizing theorems of  H. Weyl,
I. D. Berg and \linebreak[4] G. J. Murphy.

Last but not least one realizes that the
forthcomming theory is closely related to the describtion of self-dual Hilbert
AW*-modules over commutative AW*-algebras in terms of Boolean valued analysis
and logic created by  M. Ozawa and  G. Takeuti
\cite{Ozawa:83,Ozawa:85,Takeuti:78,Takeuti:83/1,Takeuti:83/2}
(\cite{Nishimura:91}) during 1979-85.
There are also relations to the work of \linebreak[4] H. Takemoto
\cite{Takemoto:73,Takemoto:75,Takemoto:76}
who has described similar phenomena in terms of continuous fields of Hilbert
spaces. In the present more special case the mathematical terminology
describing the situation is taken from measure theory.

The present paper is organized as follows: The first section is a short
summary of facts from the theory of direct integrals of measurable fields of
Hilbert spaces and of von Neumann algebras, at one side, and from the theory
of Hilbert W*-modules over commutative Hilbert W*-algebras, at the other.
We slightly modify the traditional denotations for our purposes and recall
some necessary facts from the literature. The second section deals with the
interrelation between self-dual Hilbert W*-modules over commutative
W*-algebras {\bf A} possessing a countably
generated {\bf A}-pre-dual Hilbert {\bf A}-module and special sets of mappings
into measurable fields of Hilbert spaces, giving rise to isomorphisms.
Considering von Neumann algebras of bounded module operators on those Hilbert
W*-modules we obtain their direct integral decomposition. As an application for
commutative W*-algebras {\bf A} with a special property (*) we prove that on
self-dual Hilbert {\bf A}-modules which possess
a countably generated {\bf A}-pre-dual Hilbert {\bf A}-module every normal
bounded module operator $T$ is decomposable into the sum of a normal,
diagonalizable bounded module operator $D$ and a ''compact'' bounded module
operator $K$.

The author thanks J. Friedrich, A. Kasparek, P. Richter, R. Schaflitzel and
\linebreak[4]
K. Schm\"udgen for helpful discussions and suggestions during the time of
preparation of the paper.

\newpage

\section{Preliminaries}

   We start with some necessary informations about Hilbert W*-modules.
   Throughout
   the present paper the symbol {\bf A} is denoting a C*-algebra. We make the
   convention that all modules over {\bf A} are left modules by definition.
   Following W. L. Paschke \cite{Paschke:73} and other authors
   \cite{Dupre/Gillette:83,Kasparov:80,Ozawa:83,Ozawa:85,Widom:56} we define
   a {\it pre--Hilbert {\bf A}-module} over a certain C*-algebra {\bf A} as an
   {\bf A}-module $\cal H$ equipped with a mapping $\langle .,. \rangle :
   {\cal H} \, {\rm x} \, {\cal H} \longrightarrow {\bf A}$ satisfying:

   \newcounter{raf}
   \begin{list}{(\roman{raf})}{\usecounter{raf}}
      \item $\lambda (a x) = (\lambda a) x = a (\lambda x) $ for every
      $\lambda \in {\bf C}$, $a \in {\bf A}$, $x \in \cal H$.

      \item $\langle x,x \rangle \geq 0$ for every $x \in \cal H$.

      \item $\langle x,x \rangle =0$ if and only if $x=0$.

      \item $\langle x,y \rangle = \langle y,x \rangle$* for every $x,y \in
      \cal H$.

      \item $\langle a x + b y,z \rangle = a \langle x,z \rangle + b \langle
      y,z \rangle$ for every $a,b \in {\bf A}$, every $x,y,z \in \cal H$.
   \end{list}

   The mapping $\langle .,. \rangle $ is the so called {\bf A}{\it -valued
   inner product on} $\cal H$. A pre--Hilbert {\bf A}-module is called to be
   {\it Hilbert} if it is complete with respect to the norm $\| x \| =
   \| \langle x,x \rangle \|_A^{1/2}$.
   Two Hilbert {\bf A}-modules $\{ {\cal H}_1 , \langle .,. \rangle_1 \}$,
   $\{ {\cal H}_2 , \langle .,. \rangle_2 \}$ over a certain C*-algebra {\bf A}
   are isomorphic if there exists a bijective, {\bf A}-linear, bounded mapping
   $T: {\cal H}_1 \longrightarrow {\cal H}_2$ such that $\langle .,. \rangle_1
   \equiv \langle T(.),T(.) \rangle_2$ on ${\cal H}_1 {\rm x} {\cal H}_1$.
   A Hilbert {\bf A}-module $\cal H$ is {\it finitely generated} if it is
   finitely generated as an {\bf A}-module. It is {\it countably ge\-nerated}
   if
   there exists a countable set of generators inside $\cal H$ such that the
   set of all finite {\bf A}-linear combinations of generators is norm-dense
   in $\cal H$. A Hilbert {\bf A}-module $\{ {\cal H} , \langle .,. \rangle \}$
   over a C*-algebra {\bf A} is {\it faithful} if the norm-closed
   {\bf A}-linear
   hull of the range of the inner product, $\langle \cal H, \cal H \rangle$,
   is identical with {\bf A}.

   A central notion in the theory of Hilbert C*-modules is the notion of
   self-duality since self-dual Hilbert C*-modules form a proper subcategory
   of the category of Banach C*-modules with advantageous properties, cf.
   \cite{Frank:89,Frank:91}. Denote by $\cal H'$ the set of all bounded module
   maps $f: {\cal H} \longrightarrow {\bf A}$. Following W. L. Paschke
   \cite{Paschke:73}  a Hilbert C*-module $\{ {\cal H}, \langle .,. \rangle \}$
   is called to be {\it self-dual} if every map $r \in \cal H'$ is of the form
   $\langle ., a_r \rangle$ for a certain element $a_r \in \cal H$. Let us
   remark,
   that a Hilbert AW*-module over a commutative AW*-algebra is self-dual if and
   only if it is Kaplansky-Hilbert.

   In the following we direct our attention to Hilbert W*-modules over
   commutative W*-algebras. In that case the {\bf A}-valued inner product
   on $\{ {\cal H}, \langle .,. \rangle \}$ lifts to an {\bf A}-valued inner
   product $\langle .,. \rangle_D$ on the Banach {\bf A}-module $\cal H'$
   turning $\{ {\cal H'}, \langle .,. \rangle_D \}$ into a self-dual
   Hilbert {\bf A}-module. The equalities
   \[
   \langle \langle .,x \rangle, \langle .,y \rangle \rangle_D =
   \langle x,y \rangle \, {\rm for} \, {\rm every} \, x,y \in \cal H,
   \]
   \[
   \langle \langle .,x \rangle r(.) \rangle_D = r(x) \, {\rm for} \,
   {\rm every} \, x \in {\cal H}, \, {\rm every} \, r \in \cal H'
   \]
   are satisfied, cf. \cite[Th. 3.2]{Paschke:73}. Moreover, the following
   criterion for self-duality can be formulated:

      \begin{proposition}
        [cf. {\cite[Th. 3.2, Th. 3.12]{Paschke:73}}] Let {\bf A} be a
        commutative
        W*-al\-gebra and let $\{ {\cal H}, \langle .,. \rangle \}$ be a
        Hilbert
        {\bf A}-module. Then the following two conditions are equivalent:

        \newcounter{gtr}
        \begin{list}{(\roman{gtr})}{\usecounter{gtr}}
          \item $\cal H$ is self-dual.

          \item There exist an index set $I$ and a collection of (not
          necessarily
          distinct) projections $\{ p_\alpha : \alpha \in I \}$ of {\bf A}
          indexed
          by $I$ such that $\cal H$ is isomorphic to the set of all $I$-tuples
          \[
          \tau - \Sigma \{ {\bf A}p_\alpha : \alpha \in I \} =
             \{ \{ x_\alpha \} : x_\alpha \in {\bf A}p_\alpha, \alpha \in I,
             \| \sum x_\alpha x_\alpha^* \|_A < + \infty \}
          \]
          equipped with the {\bf A}-valued inner product
          \[
          \langle x,x \rangle = w^*-\lim_{S \in \cal F} \sum_{\alpha \in S}
          x_\alpha x_\alpha^* ,
          x = \{ x_\alpha : \alpha \in I \},
          \]
          where $\cal F$ is the net of all finite subsets of $I$ being
          partially ordered by inclusion.
       \end{list}
    \end{proposition}

    \begin{corollary}
       Let {\bf A} be a commutative W*-algebra, $\cal H$ be a self-dual Hilbert
       {\bf A}-module
       being representable as $\tau - \Sigma \{ {\bf A}p_\alpha :
       \alpha \in I \}$
       for a countable set $I$. Then there exists a countably generated
       Hilbert
       {\bf A}-module $\cal K$ such that the {\bf A}-dual Banach
       {\bf A}-module
       of $\cal K$ is $\cal H$.
    \end{corollary}

    Beside the Hilbert {\bf A}-modules we would like to consider
    {\bf A}-linear
    bounded operators $T$ on them. If the underlying Hilbert C*-module
    $\cal H$ is self-dual
    they always possess an adjoint operator $T^*$ being bounded and
    {\bf A}-linear.
    The set of all such ope\-rators , ${\bf End}_A(\cal H)$, forms a
    C*-algebra in
    that situation. Moreover, ${\bf End}_A(\cal H)$ becomes a W*-algebra over
    self-dual Hilbert W*-modules, cf. \cite{Paschke:73}. An important subset
    of ${\bf End}_A(\cal H)$ is the set of {\it ''compact''} operators
    ${\bf K}_A(\cal H)$
    being defined as the norm-closed linear hull of the set
    \[
    \{ \theta_{a,b} \in {\bf End}_A({\cal H}) :\theta_{a,b}(c) = \langle c,a
    \rangle b
    \: {\rm for} \: {\rm every} \: a,b,c \in \cal H \}.
    \]
    It is a C*-subalgebra and a two-sided ideal of ${\bf End}_A(\cal H)$.

\medskip
Our standard reference sources for direct disintegration theory are the
monographs \cite{Maurin:67,Riesz:56,Takesaki:79} and the papers
\cite{Effros:66,Segal:51/1,Segal:51/2}. For recent developements in this
area see
\cite{Kehlet:78,Marechal:69,Vesterstrom/Wils:70,Wils:70},
\cite{Schaflitzel:90/1,Schaflitzel:90/2}.

\smallskip
The following definition we would like to take as a basis
      (\cite[Def. 8.9]{Takesaki:79}, \cite[p. 206-207]{Maurin:67}):
      Let $X$ be a locally compact Hausdorff measure space with Borel measure
      $\mu$. A set $\{ H_x : x \in X \}$ of Hilbert spaces indexed by $X$
      is called to be a $\mu${\it -measurable field of Hilbert spaces} if
      there exists a subspace
      $\cal E$ of the product space $\prod \{ H_x : x \in X \}$ with the
      properties:

      \newcounter{p}
      \begin{list}{(\roman{p})}{\usecounter{p}}
        \item For every $z \in \cal E$ the function $\| z(x) \|$ is an element
        of
        $L^{\infty}(X, \mu )$.

        \item If for a certain $y \in \prod \{ H_x : x \in X \}$ the function
        $\langle y(x), z(x) \rangle$ belongs to $L^{\infty}(X, \mu )$
        for every $z \in \cal E$ then $y \in \cal E$.

        \item There exists a countable subset $\{ z_i : i \in {\bf N} \}$
        of elements of $\cal E$ such that for every $x \in X$ the set
        $\{ z_i(x) : i \in {\bf N} \}$ is a basis of the Hilbert space $H_x$.
      \end{list}

   Elements $\{ h_x : x \in X \}$ of $\cal E$ are called  to be $\mu${\it
   -measurable.} We will specify the set $\cal E$ in further considerations
   since the structure of $\cal E$ is sometimes important for our purposes.
   Let us
   remark that (iii) implies the separability of the Hilbert spaces $H_x$,
   $x \in X$, of the $\mu$-measurable field $\{ H_x : x \in X \}$. Moreover,
   the map $x \in  X \longrightarrow dim(H_x) \in {\bf R} $ is
   $\mu$-measurable.

   Denote by $L^{\infty}(X, \mu , \{ H_x : x \in X \})$ the set of all rest
   classes
   of essentially bounded, $\mu$-measurable mappings of X into the
   $\mu$-measurable
   field of Hilbert spaces \linebreak[4] $\{ H_x : x \in X \}$,
   ($x \longrightarrow H_x$),
   where the elements of one rest class differ only on subsets of
   $\mu$-measure zero.
   Analogously, define  $L^1(X, \mu , \{ H_x : x \in X \})$ as the set of all
   rest
   classes of mappings $f: x \in X \longrightarrow H_x \in \{ H_x : x \in X \}$
   possessing a finite integral $\int_X \| f(x) \|_{H_x} \, d \mu (x)$ ,
   where the elements of one rest class differ only on subsets of $X$ of
   $\mu$-measure zero. Defining suitable operations, norms and other structural
   elements on \linebreak[4] $L^{\infty}(X, \mu , \{ H_x : x \in X \} )$ and
   $L^1(X, \mu , \{ H_x : x \in X \} )$
   they become a faithful self-dual Hilbert $L^{\infty}(X, \mu )$-module
   with countably generated $L^{\infty}(X, \mu )$-pre-dual Hilbert
   $L^{\infty}(X, \mu )$-module and a Banach $L^{\infty}(X, \mu )$-module,
   respectively. The third structure needed in the following is the classical
   direct integral of
   the $\mu$-measurable field of Hilbert spaces $\{ H_x : x \in X \}$,
   $\int_X H_x \, d \mu (X)$.

   Recall, that an operator $T$ on $\int_X H_x \, d \mu (x)$ is called to be
   $\mu$-{\it measurable} if the operator $T(x)$ acts on $H_x$
   as a bounded linear operator for almost every $x \in X$ and $T ({\cal E})
   \subseteq
   {\cal E}$. Now, following \cite[Def. 7.7, Cor. 7.8, Def. 7.9]{Takesaki:79}
   denote by $\int_X {\bf End}_C (H_x ) \, d \mu (x)$ the set of
   {\it decomposable operators} on $\int_X H_x \, d \mu (x)$ , i.e. the
   set of all rest classes of essentially bounded, $\mu$-measurable fields
   of operators $\{ B_x : B_x \in {\bf End}_C (H_x ); x \in X \}$, where the
   elements of one rest class differ only on subsets of X of $\mu$-measure
   zero.
   Note, that the commutant with respect to ${\bf End}_C(\int_X H_x \, d \mu
   (x))$
   of the set $\{ a \cdot id_{L^2} : a \in L^{\infty}(X, \mu ) \}$
   of all diagonal operators on $\int_X H_x \, d \mu (x))$ equals to
   $\int_X {\bf End}_C (H_x ) \, d \mu (x)$
   at least if the Hilbert spaces $H_x$ are $\mu$-almost everywhere separable.
   This property is lost in certain cases when the Hilbert spaces $H_x$ are
   taken to be non-separable, (cf. \cite{Schaflitzel:90/1,Schaflitzel:90/2}).
   With suitable chosen operations it is a normed $*$-algebra.
   Moreover, in the classical situation when $\mu$-almost all Hilbert spaces
   $H_x$ are
   separable it is a W*-algebra of type I.

    \bigskip
    Now we are prepared for further considerations.

    \newpage

\section{An equivalence principle}

Let ${\bf A}= L^{\infty}(X,\mu)$ be a commutative W*-algebra, $X$ be a
suitable
chosen locally compact, Hausdorff measure space with Borel measure $\mu$. The
purpose of the consi\-de\-rations below is , first, to show that each
self-dual Hilbert
{\bf A}-module $\cal H$ possessing a countably generated {\bf A}-pre-dual
Hilbert {\bf A}-module is isomorphic to a certain Hilbert
$L^{\infty}(X,\mu)$-module
of type $L^{\infty}(X, \mu , \{H_x : x \in X \})$ for a suitable chosen
$\mu$-measurable
field of separable Hilbert spaces $\{H_x : x \in X\}$ on $(X, \mu )$, and
secondly, to derive the direct integral decomposition of the operator algebra
${\bf End}_A(\cal H)$. Moreover, we will look for possibilities and bounds of
generalization of these equivalence relations we get.
Formulating the theorems below we enclose two results of I. E. Segal
\cite{Segal:51/1,Segal:51/2} for completeness.

   \begin{theorem}
      {\rm (existence of isomorphisms)}

      Let {\bf A} be a commutative W*-algebra and $\cal H$ be a self-dual
      Hilbert
      {\bf A}-module being the {\bf A}-dual Banach {\bf A}-module of a
      countably
      generated Hilbert {\bf A}-module.
      Then there exists a locally compact, Hausdorff measure space X with a
      Borel measure $\mu$ and a $\mu$-measurable field of Hilbert spaces
      $\{H_x : x \in X \}$
      such that:

      \newcounter{m}
      \begin{list}{(\roman{m})}{\usecounter{m}}
        \item {\bf A} is (isometricly) $*$-isomorphic to
        $L^{\infty}(X, \mu )$.
        \item $\cal H$ is isometricly isomorphic to
        $L^{\infty}(X, \mu ,\{H_x : x \in X \})$
         as a Hilbert {\bf A}-module.
        \item ${\bf End}_A(\cal H)$ is (isometricly) $*$-isomorphic to the
        W*-algebra
         \linebreak[4] $\int_X {\bf End}_C (H_x ) \, d \mu (x)$ on the Hilbert
         space
         $\int_X H_x \, d \mu (x)$.
        \item The pre-dual of $\cal H$ is isometricly isomorphic to
        $L^1(X, \mu , \{H_x : x \in X \})$.
        \item The pre-dual of ${\bf End}_A(\cal H)$ is isometricly isomorphic
        to
         \linebreak[4] $L^1(X, \mu , \{ [ {\bf End}_C(H_x) ]_* : x \in X \})$.
      \end{list}
   \end{theorem}
   \newpage

   \begin{theorem}
      {\rm (uniqueness of isomorphisms)}

      If under the assumptions of the previous theorem the W*-algebra
      {\bf A} is faithfully
      and normally representable on a separable Hilbert space then one has:

      \newcounter{nie}
      \begin{list}{(\roman{nie})}{\usecounter{nie}}
         \item If there exist two locally compact, second countable, Hausdorff
         measure spaces $X_1,X_2$ equipped with the $\sigma$-finite Borel
         measures
         $\mu_1, \mu_2$ , respectively, such that {\bf A} is $*$-isomorphic to
         both $L^{\infty}(X_1, \mu_1)$ and $L^{\infty}(X_2, \mu_2)$ then
         there exist two null sets $N_1 \subset X_1$, $N_2 \subset X_2$, a
         Borel
         isomorphism $\phi : X_2 \backslash N_2 \longrightarrow X_1 \backslash
         N_1$
         and a $*$-isomorphism $\pi : L^{\infty}(X_1, \mu_1) \longrightarrow
         L^{\infty}(X_2, \mu_2)$
         such that $\mu_1$ and $\phi (\mu_2)$ are equivalent in the sense of
         absolute continuity on $X_1 \backslash N_1$, and
         the equality $\pi (a)(x) = a(\phi (x))$ holds for every
         $a \in L^{\infty}(X_1, \mu_1 )$ and for every $x \in X_2 \backslash
         N_2$.

         \item If there are, additional, two different $\mu_{1,2}$-measurable
         fields
         of Hilbert spaces $\{ H_x^{(1)} : x \in X \}$, $\{ H_x^{(2)} :
         x \in X \}$
         satisfying condition (ii) of the previous theorem then there exist
         two null sets
         $Y_1 \subset X_1$, $Y_2 \subset X_2$ , a Borel isomorphism
         $\psi :X_2 \backslash Y_2 \longrightarrow X_1 \backslash Y_1$ and a
         $\mu_1$-$\mu_2$-measurable field of unitary operators
         $\{ U_x : H_x^{(1)} \longrightarrow H_{\psi^{-1} (x)}^{(2)} : x \in
         X_1 \backslash Y_1 \}$
         such that $\mu_1$ and $\psi (\mu_2)$ are equivalent in the sense of
         absolute
         continuity on $X_1 \backslash Y_1$, and that $U_x {\bf End}_C
         (H_x^{(1)}) U_x^*=
         {\bf End}_C(H_{\psi^{-1} (x)}^{(2)})$ for every $x \in X_1
         \backslash Y_1$.
      \end{list}
   \end{theorem}

{\it Proofs of the theorems}: The assertion (i) of the first theorem was
proved
by I.E.Segal (\cite{Segal:51/1,Segal:51/2}) in the early fiftees and can
be found at \cite[Th. 3.4.4]{Riesz:56}, whereas item (i) of the second theorem
can be derived from \cite[Lemma 8.22, Th. 8.23]{Takesaki:79} as a special case.
Therefore, one can identify {\bf A} with $L^{\infty}(X_K, \mu_K )$ for a
special
locally compact, Hausdorff measure space $X_K$ with Borel measure $\mu_K$
being constructed from the compact Hausdorff space $K$ realizing the
$*$-isomorphy
${\bf A}=C(K)$ along the line of \cite[p.110]{Takesaki:79};
i.e., taking $X_K$ as the union of the support sets $\Gamma_\alpha \subseteq
K$
of a maximal family of positive normal measures $\mu_\alpha$ on $K$ with
disjoint supports, and defining $\mu_K$ on $X_K$ by the formula $\mu_K(f) =
\sum_\alpha \mu_\alpha (f)$ for each continuous function $f$ on $X_K$ with
compact support. Finally, one has a bijection between continuous functions on
$K$ and rest classes of $\mu_K$-measurable, essentially bounded functions of
$L^{\infty}(X_K, \mu_K)$. Now, according to the isometric isomorphy of the
Hilbert {\bf A}-modules ${\cal H} = \tau - \Sigma \{ {\bf A}p_i : i \in
{\bf N} \}$
one gets assertion (ii) of the first theorem identifying {\bf A} with $C(K)$
and considering the Hilbert spaces $H_x = \tau-\Sigma \{ f(x) \cdot p_i(x) :
f \in C(K), i \in {\bf N} \}$ for each $x \in X_K \subseteq K$. They form a
$\mu_K$-measurable
field of separable Hilbert spaces on $X_K$, where $\cal E$ is the subset of
all
square-integrable elements of $L^\infty(X_K, \mu_K, \{ H_x : x \in X_K \} )$.
Note, that $\cal E$ is norm-dense in $\int_{X_K} H_x \; d\mu_K(x)$ by
definition
and $\tau_1$-dense in $L^\infty(X_K, \mu_K, \{ H_x : x \in X_K \} )$.
Turning to ${\bf A}= L^{\infty}(X_K, \mu_K)$
and recalling the definition of $\tau - \Sigma$ type Hilbert C*-modules
one finishes.

Now consider the set ${\bf End}_A(\cal H)$ of all bounded, {\bf A}-linear
operators on $\cal H$. The boundedness and the {\bf A}-linearity of these
operators
guarantee the invariance of the Hilbert spaces $H_x$ under the action of them.
Moreover, ${\bf End}_A(\cal H)$ acts on each Hilbert space $H_x$,
($x \in X_K$),
like ${\bf End}_C(H_x)$ and, globally, preserves $\cal E$ and the
$\mu_K$-measurability
of the field
$\{ H_x : x \in X_K \}$. That is, ${\bf End}_A(\cal H)$ is embeddable into
$\int_X {\bf End}_C (H_x ) \, d \mu (x)$ in the sense of the coincidence of
these operators on $\cal E$. Vice versa, every
essentially bounded, $\mu_K$-measurable map from $X_K$ onto $\{ {\bf End}_C
(H_x) : x \in X_K \})$
induces a bounded, {\bf A}-linear operator on $\cal E$ and, hence, on
$\cal H$. So one has shown the
(isometric) $*$-isomorphy of these W*-algebras, i.e. assertion (iii) of the
first theorem.

Item (ii) of the second theorem can be derived from the isometric Hilbert
{\bf A}-module isomorphism of $L^{\infty}(X_1, \mu_1, \{ H_x^{(1)} :
x \in X_1 \})$
and $L^{\infty}(X_2, \mu_2, \{ H_x^{(2)} : x \in X_2 \})$, from the
commutativity of {\bf A} and from (i) of both the theorems, whereas
the facts
(iv) and (v) of the first theorem follow from the self-duality of Hilbert
spaces and from \cite[p.70, Prop.]{Riesz:56} or from
\cite[Prop. 8.38]{Takesaki:79} ,
respectively.  $\bullet$

\begin{example}    {\rm
   For a fixed W*-algebra {\bf A} use the denotations
   \[
   l_2({\bf A}) = \left\{ \{a_i\}_{i \in {\bf N}} : a_i \in {\bf A},
                  {\sum}_i a_ia_i^* \, {\rm converges} \, {\rm  in} \,
                  \|.\|_A             \right\}
   \]
   \[
   l_2({\bf A}){\rm '} = \left\{ \{a_i\}_{i \in {\bf N}} : a_i \in {\bf A},
                  \left\| {\sum}_i a_ia_i^* \right\|_A \, {\rm converges} \,
                        \right\}
   \]
   for the standard countably generated Hilbert {\bf A}-module and its
   {\bf A}-dual Banach {\bf A}-module.

      \newcounter{po}
      \begin{list}{\alph{po})}{\usecounter{po}}
         \item Let ${\bf A}= l^{\infty}$ and ${\cal H} = (l^2(l^{\infty}))'$.
         The W*-algebra $l^{\infty}$ is faithfully representable as the von
         Neumann algebra of all bounded, diagonal operators on the separable
         Hilbert space $l^2$. According to the theorems above one has
         \[ {\cal H} = (l^2(l^{\infty}))' = L^{\infty}({\bf N}, \nu,
         \{ l^2_{(i)} : i \in {\bf N} \}), \]
         \[ {\bf End}_A({\cal H}) = \int_{\bf N} {\bf End}_C (H_x ) \, d \nu
         (x), \]
         where $\nu$ denotes a discrete measure on the set of natural numbers
         {\bf N}.

         \item For ${\bf A} = L^{\infty}([0,1] , \lambda)$ and
         ${\cal H}=(l^2(L^{\infty}([0,1] , \lambda)))'$ one has
         \[
         {\cal H} = L^{\infty}([0,1] , \lambda, \{ l^2_{(i)} : i \in [0,1]  \})
         \]
         \[
         {\bf End}_A({\cal H}) = \int_{[0,1]} {\bf End}_C (H_x ) \, d \lambda
         (x),
         \]
         where $\lambda$ denotes the Lebesgue measure on the unit interval.

      \end{list}                      }
\end{example}

   \begin{corollary}[cf. {\cite[Cor. 8.20]{Takesaki:79}}]
      For ${\bf A} = L^{\infty}(X, \mu )$,
      ${\cal H}= (l^2({\bf A}))'$ one has \linebreak[4] ${\bf End}_A({\cal H})
      = {\bf End}_C(l^2)
      \overline{\otimes} {\bf A}$, where $\overline{\otimes}$
      denotes the w*-tensor product.
   \end{corollary}

   \begin{remark}
      {\rm If the locally compact measure space $X$ is not second countable and
      the
      Borel measure $\mu$ on $X$ is not $\sigma$-finite then the statement of
      Theorem 2.2  is not longer true, in general. If one omits the
      separability
      condition to the Hilbert spaces $H_x$ then Theorem 2.1 fails to be true,
      in general.

      For example, consider the W*-algebra ${\bf A}=L^{\infty}([0,1], \lambda
)$,
      where $\lambda$ denotes the Lebesgue measure, and the Hilbert
      {\bf A}-module
      ${\cal H}=\tau - \Sigma \{ {\bf A}_{(\alpha)} : \alpha \in \{ [0,1],
      \nu \} \}$,
      with $\nu$ a discrete measure on $[0,1]$. Following the idea of Theorem
      2.1,(ii)
      one should like to compare $\cal H$ with $L^{\infty}([0,1], \lambda ,
      \{ l^2_{(\alpha)}([0,1]) : \alpha \in \{ [0,1], \nu \} \})$.
      (For the more complicated general definition of a $\mu$-measurable field
      of non-separable Hilbert spaces see R. Schaflitzel
      \cite{Schaflitzel:90/1,Schaflitzel:90/2} e.g..) However, the special
      mapping
      \begin{eqnarray*}
      \{ [0,1], \lambda \} & \longrightarrow & l^2([0,1]) \\
      x & \longrightarrow & \{ \delta_{x,t}(t) : t \in \{ [0,1], \nu \} \}
      \end{eqnarray*}
      belongs to $L^{\infty}([0,1], \lambda , \{ l^2_{(\alpha)}([0,1]) :
      \alpha \in \{ [0,1], \mu \} \} )$
      as a non-zero element, whereas its reflection inside ${\cal H}=
      \tau - \Sigma \{ {\bf A}_{(\alpha)} : \alpha \in \{ [0,1], \nu \} \}$
      gives the zero element.

      Beside this,  from a result
      of R. Schaflitzel \cite{Schaflitzel:90/1},
      \cite[Lemma 6]{Schaflitzel:90/2}
      there follows that
      under the assumption of the continuum-hypothesis it may happen that the
      algebra of decomposable operators is not the
      commutant of the algebra of diagonalizable operators on  direct integrals
      of
      certain $\mu$-measurable fields of Hilbert spaces \linebreak[4]
      $\{ H_x : x \in X \}$ on $X$ with dim$(H_x) \geq {\bf c}$,
      $card(X) \geq {\bf c}$ and an almost pointwise ortho\-gonal generating
      set
      $\Gamma_o$ of the corresponding direct integral of the Hilbert spaces
      $H_x$ with card$(\Gamma_o) > {\bf c}$. As a concrete example he
      considered
      $X = [0,1]$, $\lambda $ - the Lebesgue measure and $H_x = l^2([0,1])$,
      i.e.
      the same situation as above. }
   \end{remark}

Now we are interested in a direct integral decomposition of von Neumann
algebras {\bf M} of operators on self-dual Hilbert W*-modules $\cal H$ over
commutative W*-algebras
{\bf A}. Conversely, we ask for which W*-subalgebras {\bf B} of the centre of
a given W*-algebra {\bf M} there exists a self-dual Hilbert {\bf B}-module
$\cal H$ such that {\bf M} is faithfully $*$-representable as a von Neumann
subalgebra of ${\bf End}_B(\cal H)$. The answer can be derived from the
direct
disintegration theory, especially from \cite[Th. 8.22, Th. 8.23]{Takesaki:79}.

   \begin{corollary}
      Let {\bf A} be a commutative W*-algebra and $\cal H$ be a self-dual
      Hilbert
      {\bf A}-module possessing a countably generated {\bf A}-pre-dual Hilbert
      {\bf A}-module. Let \linebreak[4] ${\bf M} \in {\bf End}_A(\cal H)$ be a
      von Neumann subalgebra.

      If ${\bf A}=L^{\infty}(X, \mu)$, ${\cal H}=L^{\infty}(X, \mu ,
      \{ H_x : x \in X \})$
      in the sense of Theorem 2.1,(i),(ii) then there exists a
      $\mu$-measurable field
      of von Neumann algebras \linebreak[4] $\{ {\bf M}_x : x \in X \}$ on
      the $\mu$-measurable
      field of Hilbert spaces $\{ H_x : x \in X \}$ such that {\bf M} is
      (isometricly) $*$-isomorphic to $\int_X {\bf M}_x \, d \mu (x)$.
   \end{corollary}

   \begin{corollary}
      Let {\bf M} be a W*-algebra possessing a normal, faithful representation
      on a separable Hilbert space. Let {\bf B} be a W*-subalgebra of the
      centre of {\bf M} and let {\bf B} be $*$-isomorphic to $L^{\infty}(X,
      \mu )$
      for a certain locally compact, Hausdorff, second countable measure
      space $X$ with $\sigma$-finite Borel measure $\mu$. Then for
      $\mu$-almost
      every $x \in X$ there exist a Hilbert subspace $H_x \subseteq H$ and a
      von Neumann algebra ${\bf M}_x \subseteq {\bf End}_C(H_x)$ such that:

      \newcounter{r}
      \begin{list}{(\roman{r})}{\usecounter{r}}
      \item The set of Hilbert spaces $\{ H_x : x \in X \}$ is a
      $\mu$-measurable
      field of Hilbert spaces, and $\int_X H_x \, d \mu (x) = H$.

      \item The set of von Neumann algebras $\{ {\bf M}_x : x \in X \}$ is a
      $\mu$-measurable field and {\bf M} is $*$-isomorphic to
      $\int_X {\bf M}_x \, d \mu (x)$.

      \item {\bf M} is faithfully representable as a von Neumann algebra of
      bounded module operators on the self-dual Hilbert {\bf B}-module
      ${\cal H}=L^{\infty}(X, \mu , \{ H_x : x \in X \})$ being the
      {\bf B}-dual
      Banach {\bf B}-module of a countably generated Hilbert {\bf B}-module.
      \end{list}

   \end{corollary}

Analysing these statements one concludes that the theory of direct integral
decomposition of W*-algebras possessing a normal, faithful representation on a
separable Hilbert space is one-to-one translatable to the theory of von Neumann
subalgebras of the W*-algebras of bounded module operators on self-dual Hilbert
{\bf A}-modules with countably generated {\bf A}-pre-dual Hilbert
{\bf A}-module
over separable commutative W*-algebras {\bf A}.

\section{An application}

 In the present section we like to generalize the following theorem of
 H. Weyl and \linebreak[4] I.D. Berg (\cite{Weyl:09} 1909 and \cite{Berg:71}
 1971):

    \begin{proposition}
       [{\rm  Berg's and  Weyl's theorem}]
       Every linear bounded normal ope\-rator $T$ on a separable Hilbert space
       is
       decomposable into the sum of a normal,diago\-na\-lizable and a compact
       operator, $D$ and $K$.
       If $T$ is self-adjoint then for every given \linebreak[4] $\varepsilon
       > 0$ one can even choose
       self-adjoint $D$ and $K$ such that $\|K\| < \varepsilon$.
    \end {proposition}

 In 1970  P. R. Halmos has shown by some examples that  Weyl's theorem can not
 be generalized for self-adjoint operators on non-separable Hilbert spaces,
 (\cite{Halmos:70}).
 Nevertheless, the two theorems of the previous paragraph and a result of
 \linebreak[4]
 G. J. Murphy \cite[Th. 9]{Murphy:88} suggest to us another
 way of generalization weakening the notion of compactness.

 Let {\bf A} be a commutative W*-algebra and $\cal H$ be a self-dual Hilbert
 {\bf A}-module with countably generated {\bf A}-pre-dual Hilbert
 {\bf A}-module.
 Let us call an operator $T \in {\bf End}_A(\cal H)$ to be {\it diagonalizable}
 if
 and only if there exist a sequence of pairwise orthogonal projections
 $\{ P_i : i \in {\bf N} \}$ of ${\bf K}_A(\cal H)$ and a sequence of elements
 $\{ a_i : i \in {\bf N} \}$
 of {\bf A} such that $T = \sum_{i \in  N} a_i P_i$ in the sense of
 w*-convergence.
 Furthermore, we say that the commutative W*-algebra {\bf A} has
 {\it property (*)}
 if and only if the set of all normal states $f$ on {\bf A} with range
 projection $p_f$,
 for which the norm completion of the pre--Hilbert space
 $\{ {\bf A}p_f, f(\langle .,. \rangle_A) \}$
 is separable, separates the elements of {\bf A}. Finally, we call a locally
 compact Hausdorff measure space $X$ to be {\it locally second countable} if
 for
 every $x \in X$ there exists a clopen subset $Y \subseteq X$ containing $x$
 and being second countable with respect to the measure $\mu$.
 We get the following result using assertions of R. V. Kadison
 \cite{Kadison:83,Kadison:84}
 and of K. Grove, G. K. Pedersen \cite{Grove/Pedersen:84} in the proof:

    \begin{theorem}
       Let {\bf A} be a commutative W*-algebra with property (*). Let $\cal H$
       be a self-dual Hilbert {\bf A}-module possessing a countably generated
       {\bf A}-pre-dual Hilbert {\bf A}-module.
       Let $T$ be an {\bf A}-linear bounded normal operator  on $\cal H$.

       \noindent
       Then $T$ is decomposable into the sum
       of a {\bf A}-linear bounded normal diagonalizable operator D on $\cal H$
       and a {\bf A}-linear bounded ''compact'' operator K on $\cal H$.
       If $T$ is self-adjoint then for every $\varepsilon > 0$ the operators
       $D$ and $K$
       can be chosen to be self-adjoint and such that $\|K\| < \varepsilon$.
    \end{theorem}

    \begin{corollary}
        Let $X$ be a locally compact, locally second countable Hausdorff
        measure space with Borel measure $\mu$. Let $\{ H_x : x \in X \}$ be
        a $\mu$-measurable field of Hilbert spaces on $X$. Then each normal
        decomposable bounded linear operator $T$ on the (non-separable, in
        general)
        Hilbert space $H=L^2(x, \mu , \{ H_x : x \in X \})$ can be decomposed
        into
        the sum of a normal diagonalizable decomposable bounded linear
        operator $D$
        on $H$ and a decomposable bounded linear operator $K$ on $H$,
        $K$ being compact on every subspace $L^2(Y, \mu , \{ H_x : x \in Y \})
        \subseteq H$ with $Y \subseteq X$ being clopen and second countable.
    \end{corollary}

{\it Proof:} Choose a representation of the commutative W*-algebra {\bf A} as
$L^{\infty}(X, \mu )$ for a certain locally compact Hausdorff measure space
$X$
with Borel measure $\mu$. By the assumptions $X$ is locally second countable.
Since $X$ is the union of a family of second countable clopen subsets
$Y_\alpha$
with pairwise empty intersection one can suppose without loss of  generality
that $X$ is second countable and, consequently, $\mu$ is $\sigma$-finite.

By the first theorem of the previous section there exists a $\mu$-measurable
field
of Hilbert spaces $\{ H_x : x \in X \}$ such that $\cal H$ is isometricly
isomorphic to \newline
$L^{\infty}(X, \mu , \{ H_x : x \in X \})$ and that ${\bf End}_A(\cal H)$ is
$*$-isomorphic to $\int_X {\bf End}_C (H_x ) \, d \mu (x)$.
The latter can be interpreted as a C*-subalgebra of the set
of all bounded linear operators {${\bf End}_C(H)$} on the separable Hilbert
space
$H=\int_X H_x \, d \mu (x)$ by (iii) of that theorem. Consequently,
one can apply Murphy's theorem (\cite[Theorem 9]{Murphy:88}) to $T \in
{\bf End}_A({\cal H}) \equiv \int_X {\bf End}_C (H_x ) \, d \mu (x)$ under
that
point of view, and one
gets a diagonalizable normal operator $D \in \int_X {\bf End}_C (H_x ) \, d
\mu (x)
\equiv {\bf End}_A(\cal H)$ and a compact on $H$ operator $K \in
\int_X {\bf K}_C (H_x ) \, d \mu (x) \equiv {\bf K}_A(\cal H)$ such that
$T=D+K$. Pay attention, that the diagonalizability of
$D \in \int_X {\bf End}_C (H_x ) \, d \mu (x)$ on $H$ means that there exist a
countably
number of complex eigen-values $\{ \lambda_n : n \in {\bf N} \}$ and a
countably
number of projections $\{ P_n \} \in \int_X {\bf End}_C (H_x ) \, d \mu (x)
\equiv
{\bf End}_A({\cal H})$ such that $D = \sum_n \lambda_n P_n$ on $H$ and on
$\cal H$
simultaneously, i.e. $D$ is diagonalizable on $\cal H$, too. The nature of the
eigen-vectors does not matter. Moreover, if
$T$ is self-adjoint then for every $\varepsilon > 0$ one can choose $D$ and
$K$ in such a way that they are self-adjoint and  $\|K\| < \varepsilon$.

Remark, that if for certain projections $p \in {\bf A}$
the Hilbert $p{\bf A}$-module $p \cal H$ is finitely generated (or,
equivalently, $p {\bf End}_A(\cal H)$$=p{\bf K}_A(\cal H))$ then $pK=0$ and
$pT$ is diagonalizable by the results of R. V. Kadison and K. Grove, G. K.
Pedersen cited above.  $\bullet$

The corollary can be derived from the theorem using item
(iii) of Theorem 2.1 and Murphy's theorem.

   \begin{remark}
      {\rm Unfortunately, we are not able to say anything about the
      possibly validity of the theorem without assuming {\bf A}
      to have property (*).

      Beside this, a generalization to the case
      of {\bf A} being a commutative AW*-algebra with a similar property like
      property (*) in the W*-case seems to be possible using e. g. a transfer
      principle
      developed by G. Takeuti and M. Ozawa
      \cite{Ozawa:83,Ozawa:85,Takeuti:78,Takeuti:83/1,Takeuti:83/2},
      \cite{,Nishimura:91} between
      the theory of self-dual Hilbert AW*-modules over commutative AW*-algebras
      and its description in terms of Boolean valued analysis and logic.
      However, in the light of results of K. Grove and G. K. Pedersen
      \cite{Grove/Pedersen:84} much more general commutative
      C*-algebras than arbitrary AW*-algebras can not appear.
      A result of R. V. Kadison
      \cite{Kadison:83,Kadison:84} who proved that each normal element of
      the W*-algebra ${\bf M}_n({\bf A}) = {\bf End}_A({\bf A}^n)$ with
      {\bf A} being a W*-algebra is diagonalizable for every natural number $n$
      encourages to check the non-commutative case. But all that remains for
      further research.

      \medskip
      After this paper has circulated as a preprint the author had fruitful
      discussions with R. Schaflitzel about possibilities of application of the
      obtained equivalence principle to get an alternative definition of
      generalized direct integrals (i.e., the non-separable case).
      Let $I$ be an index set of non-countable cardinality card($I$). The
      self-dual Hilbert $L^{\infty}(X,\mu )$-module
      ${\cal M}_{{\rm card(}I{\rm )}} = \tau - \Sigma \{ L^{\infty}(X, \mu
      )_{(\alpha)} :
      \alpha \in I \}$ is card($I$)-homogenous for each cardinality card($I$).
      Moreover, the cardinality card($I$) of the generating set $I$ of
      ${\cal M}_{{\rm card(}I{\rm )}}$ is uniquely defined up to isomorphy of
      Hilbert C*-modules, (cf. \cite[\S 10, Th. 4]{Kaplansky}).
      Now, the principal idea is to use the existing isomorphy
      ${\cal M}_{{\rm card(}I{\rm )}} \cong L^{\infty}(X, \mu , \{ H_x : x
      \in X \})$
      (where $\{ H_x : x \in X \}$ is a certain $\mu$-measurable field of
      Hilbert spaces $H_x$ of dimension card($I$)) to define a standard direct
      integral of non-separable Hilbert spaces. Simply, take the subset of all
      square-integrable elements of $L^{\infty}(X, \mu , \{ H_x : x \in X \})$
      and close it up with respect to the direct integral norm. What turns out?
      One gets the smallest (non-separable, by construction) Hilbert space $H$
      satisfying a generalized definition for direct integrals (cf. \S 1) and
      containing the constant mappings $x \in X \rightarrow h=const. \in H_x$.
      That is , $H$ equals to the direct integral norm closure of the set of
      such elements $h \in \prod_{x \in X} H_x$ satisfying two properties:

      \smallskip \noindent
      $\;$ (i) The mapping $x \in X \rightarrow \| h(x) \|^2$ is integrable.

      \noindent
      $\;$ (ii) There is a subset $N \subseteq X$ of $\mu$-measure zero such
      that
      the set $\{h(x) : x \in X \}$ \linebreak[4] $\;$ generates a separable
      subspace of $H_x$.

      \smallskip \noindent
      Of course, since $H$ is the smallest direct integral in a certain
      sense $H$ is unique up to isomorphy.

      Consequently, one could not expect to get much more information about
      genera\-lized direct integrals in the non-separable case using only the
      described equivalence.

      }
  \end{remark}

\medskip
Universit\"at Leipzig

FB Mathematik/Informatik

Mathematisches Institut

Augustusplatz 10

D--04109 Leipzig

Fed. Rep. Germany.

frank@mathematik.uni-leipzig.d400.de


\begin{thebibliography}{99}
   \bibitem{Berg:71} I. D. Berg, An extension of Weyl--von Neumann theorem
   to normal
   operators. Trans. Amer. Math. Soc. {\bf 160}(1971), 365-371.
   \bibitem{Dupre/Gillette:83} M. J. Dupr\'e, R. M. Gillette, Banach bundles,
   Banach modules
   and automorphisms of C*-algebras. Boston--London--Melbourne: Pitman,
   Pitman Adv. Publ. Program 1983.
   \bibitem{Effros:66} E. G. Effros, Global structure in von Neumann algebras.
   Trans. Amer. Math. Soc. {\bf 121}(1966), 434-454.

   \bibitem{Frank:89} M. Frank, Self-duality and C*-reflexivity of Hilbert
   C*-modules.
   Zeitschr. Anal. Anw. {\bf 9}(1990), 165-176.
   \bibitem{Frank:91} M. Frank, Hilbert C*-modules over monotone complete
   C*-algebras
   and a Weyl--Berg type theorem. NTZ-preprint 3/91, Universit\"at Leipzig,
   FRG, 1991.
   \bibitem{Grove/Pedersen:84} K. Grove, G. K. Pedersen, Diagonalizing
   matrices over C(X).
   J. Funct. Anal. {\bf 59}(1984), 64-89.

   \bibitem{Halmos:70} P. R. Halmos, Ten problems in Hilbert spaces. Bull.
   Amer. Math. Soc. {\bf 76}(1970), 887-933.

   \bibitem{Henrichs:74} R. W. Henrichs, Maximale Integralzerlegungen
   invarianter
   positiv definiter Funktionen auf diskreten Gruppen. Math. Ann.
   {\bf 208}(1974),
   15-31.
   \bibitem{Henrichs:79} R. W. Henrichs, On decomposition theory for unitary
   representations of locally compact groups. J. Funct. Anal. {\bf 31}(1979),
   101-114.
   \bibitem{Henrichs:82} R. W. Henrichs, Decomposition of invariant states
   and
   nonseparable C*-algebras. Publ. Res. Inst. Math. Sci. {\bf 18}(1982),
   159-181.

   \bibitem{Kadison:83} R. V. Kadison, Diagonalizing matrices over operator
   algebras. Bull. Amer.
   Math. Soc. {\bf 8}(1983), 84-86.
   \bibitem{Kadison:84} R. V. Kadison, Diagonalizing matrices. Amer. J. Math.
   {\bf 106}(1984), 1451-1468.
   \bibitem{Kaplansky} I. Kaplansky, Algebras of type I. Ann. Math.
   {\bf 56}(1952),
   460-472.

   \bibitem{Kasparov:80} G. G. Kasparov, Hilbert C*-modules: Theorems of
   Stinespring and
   Voiculescu. J. Oper. Theory {\bf 4}(1980), 133-150.
   \bibitem{Kehlet:78} E. T. Kehlet, Disintegration theory on a constant
   field of
   non-separable Hilbert spaces. Math. Scand. {\bf 43}(1978), 353-362.
   \bibitem{Marechal:69} O. Mar\'echal, Champs mesurables d'espaces
   hilbertiens.
   Bull. Sc. Math. {\bf 93}(1969), 113-143.

   \bibitem{Maurin:67} K. Maurin, Methods in Hilbert space. Warszawa, Polish
   Scientific Publishers, 1967.
   \bibitem{Murphy:88} G. J. Murphy, Diagonality in C*-algebras. Math.
   Zeitschr. {\bf 199}(1988),
   279-284.
   \bibitem{Nishimura:91} H. Nishimura, Some connections between Boolean
   valued analysis and topological
   reduction theory for C*-algebras. Zeitschr. Math. Logik Grundlagen Math.
   {\bf 36}(1990), 471-479.

   \bibitem{Ozawa:83} M. Ozawa, Boolean valued interpretation of Hilbert space
   theory. J. Math. Soc. Japan {\bf 35}(1983), 609-627.
   \bibitem{Ozawa:85} M. Ozawa, Boolean valued interpretation of Banach space
   theory and module structures of von Neumann algebras. Nagoya Math. J.
   {\bf 117}(1990),
   1-36.
   \bibitem{Paschke:73} W. L. Paschke, Inner product modules over B*-algebras.
   Trans. Amer. Math. Soc. {\bf 182}(1973), 443-468.

   \bibitem{Richter:91} P. Richter, On the separability of Hilbert spaces.
   Wiss. Z. Univ. Leipzig, Math.--Nat.wiss. Reihe {\bf 39}(1990), 666-669.

   \bibitem{Riesz:56} F. Riesz, B. Sz.--Nagy, Vorlesungen \"uber
   Funktionalanalysis. Berlin:
   Verlag der Wissenschaften 1956.
   \bibitem{Schaflitzel:90/1} R. Schaflitzel, Direct integrals of not
   necessarily separable Hilbert spaces , in: Abstracts, Short communications,
   ICM Kyoto, Japan, August 21-29, 1990, p. 133.
   \bibitem{Schaflitzel:90/2} R. Schaflitzel, The algebra of decomposable
   operators
   in direct integrals of not necessarily separable Hilbert spaces. Proc.
   Amer. Math. Soc.
   {\bf 100}(1990), 983-987.

   \bibitem{Segal:51/1} I. E. Segal, Equivalence of measure spaces. Amer. J.
   Math.
   {\bf 73}(1951), 275-313.
   \bibitem{Segal:51/2} I. E. Segal, Decomposition of operator algebras. Mem.
   Amer. Math. Soc. {\bf 9}(1951), no. 1, 1-67, no.2, 1-66.
   \bibitem{Takemoto:73} H. Takemoto, On a characterization of AW*-modules and
   a representation of Gelfand type of noncommutative operator algebras,
   Michigan Math. J. {\bf 20}(1973), 115-127.
   \bibitem{Takemoto:75} H. Takemoto, Decomposable operators in continuous
   fields of Hilbert spaces.
   T{\^o}hoku Math. J. {\bf 27}(1975), 413-435.
   \bibitem{Takemoto:76} H. Takemoto, On the weakly continuous constant field
   of Hilbert
   space and its application to the reduction theory of von Neumann algebra,
   T{\^o}hoku Math. J. {\bf 28}(1976), 479-496.

   \bibitem{Takesaki:79} M. Takesaki, Theory of operator algebras,I.
   New York--Heidelberg--Berlin:
   Springer--Verlag 1979.
   \bibitem{Takeuti:78} G. Takeuti, Two applications of logic to mathematics.
   Tokyo--Princeton: Iwanami and Princeton University Press 1978.
   \bibitem{Takeuti:83/1} G. Takeuti, C*-algebras and Boolean valued analysis.
   Jap. J. Math. {\bf 9}(1983), 207-246.

   \bibitem{Takeuti:83/2} G. Takeuti, Von Neumann algebras and Boolean valued
   analysis.
   J. Math. Soc. Japan {\bf 35}(1983), 1-21.
   \bibitem{Valette:82} A. Valette, Extensions of C*-algebras: A survey of
   the
   Brown--Douglas--Fillmore theory. Nieuw Archief voor Wiskunde(3)
   {\bf 30}(1982),
   41-69.
   \bibitem{Vesterstrom/Wils:70} J. Vesterstrom, W. Wils, Direct integrals of
   Hilbert spaces II. Math. Scand. {\bf 26}(1970), 89-102.

   \bibitem{Weyl:09} H. Weyl, \"Uber beschr\"ankte quadratische Formen, deren
   Differenz vollstetig ist.
   Rend. Circ. Mat. Palermo {\bf 27}(1909), 373-392.
   \bibitem{Widom:56} H. Widom, Embedding in algebras of type I. Duke Math. J.
   {\bf 23}(1956), 309-324.
   \bibitem{Wils:70} W. Wils, Direct integrals of Hilbert spaces I., Math.
   Scand. {\bf 26}(1970),
   73-88.

\end{thebibliography}
\end{document}